\documentclass[french,manmat,final]{svjour}
\usepackage[english]{babel}
\usepackage{latexsym}
\usepackage{graphicx}
\usepackage{epsfig}
\newcommand{\proofend}{\hfill $\square$}

\usepackage{amsmath}
\usepackage{amsfonts}
\usepackage{amssymb}
\usepackage{verbatim}
\usepackage{xy}

\newcommand{\Pn}{\mathbb{P}^n}
\newcommand{\An}{\K^n}
\newcommand{\AutKn}{Aut(n,\K)}

\newcommand{\GLn}{GL(n,\K)}
\newcommand{\AffKn}{A\!f\!f(n,\K)}
\newcommand{\DKn}{D(n,\K)}
\newcommand{\ADKn}{AD(n,\K)}

\newcommand{\PGLnpl}{PGL(n+1,\K)}
\newcommand{\BirPn}{Bir(\Pn)}
\newcommand{\BirAn}{Bir(\An)}
\newcommand{\K}{\mathbb{K}}
\newcommand{\GLnZ}{GL(n,\mathbb{Z})}
\newcommand{\GLnmZ}{GL(n-1,\mathbb{Z})}
\newcommand{\hide}[1]{}

\newcommand{\KerL}[1]{\Delta_{#1}}
\newcommand{\Mt}{\ ^{t}\! M}

\newcommand{\leftexp}[2]{{\vphantom{#2}}^{#1}{\! \! \! #2}}

\newcommand{\defn}[1]{\textit{#1}}

\titlerunning{Affine automorphisms of $\An$ and linear automorphisms of $\Pn$}
\title{{\bf {\Large Conjugacy classes of affine automorphisms of $\K^n$ and linear automorphisms of $\Pn$ in the Cremona groups}} }
\author{J\'er\'emy Blanc}
\institute{J\'er\'emy Blanc \at 
Universit\' e de Gen\`eve,
Section de math\'ematiques, 
2-4 rue du Li\`evre,
CP 64, 1211 Gen\`eve 4 (Switzerland), \email{Jeremy.Blanc@math.unige.ch}}
\date{2005}
\begin{document}
\maketitle
\begin{abstract}We describe the conjugacy classes of affine automorphisms in the group $Aut(n,\K)$ (respectively $Bir(\K^n)$) of automorphisms (respectively of birational maps) of $\K^n$. From this we deduce also the classification of conjugacy classes of automorphisms of $\Pn$ in the Cremona group $Bir(\K^n)$.\end{abstract}

\section{Introduction}
Let $\K$ be an algebraically closed field of characteristic $0$ and let $\An$ and $\Pn$ denote respectively the affine and projective $n$-spaces over $\K$.

The {\it Cremona group}, which is the group of birational maps of these two spaces, $\BirAn=\BirPn$, has been studied a~lot, especially in dimension $2$ and $3$, see for example \cite{bib:HH} %\cite{bib:SC}
  and \cite{bib:AC}.
Its subgroup of biregular morphisms (or automorphisms) of $\An$, called the {\it affine Cremona group} $\AutKn$, has been also  much explored. We refer to \cite{bib:HK} for a list of references.

In both cases, the question of the conjugacy classes of elements is natural. For $Bir(\mathbb{P}^2)$, the classical approach can be found in \cite{bib:SK} and \cite{bib:AW}. A modern classification of birational morphisms of prime order was completed in \cite{bib:BaB}, \cite{bib:DF} and \cite{bib:BBl}.
We refer to \cite{bib:KS} and their references for the group $\AutKn$.

In the literature, the affine and projective cases are often treated separately, with different methods, although the groups are very close as we can see in the following diagram:
\begin{center}
$\begin{array}{ccccccccc}
 & & \PGLnpl &  & \subset &  & \BirPn \\
  & & \cup  & & & & \| \\
 \GLn &\subset &\AffKn& \subset &\AutKn& \subset &\BirAn.
\end{array}$
\end{center}
In this paper, we restrict ourselves neither to small dimensions nor to finite elements, but to the case of maps of degree $1$. We give the conjugacy classes, in the Cremona groups, of automorphisms of $\An$ and $\Pn$ that maps lines to lines. Explicitely these are the group $\AffKn$ of affine automorphisms of $\An$ and the group $\PGLnpl$ of automorphisms of $\Pn$ (or linear birational maps). These two groups are very classical and studied in many domains of mathematics.

More precisely,\textit{ the goal of this paper is to find under which conditions two affine automorphisms of $\An$ are conjugate in $\AffKn$ itself, in $\AutKn$, or in $\BirAn$.}
For this purpose, we describe the trace on $\AffKn$ of the conjugacy classes of an affine automorphism in $\AffKn$, $\AutKn$ and $\BirAn$, respectively in Sections $\ref{Sec:ConjugationAffAn}$, $\ref{Sec:ConjugationAutKn}$, and $\ref{Sec:ConjugationBirAn}$. 

The similar question in the projective case is to find under which conditions two linear automorphisms of $\Pn$ are conjugate in the group $\BirPn$. This will be answered in Section $\ref{Sec:ConjBirPn}$.

\bigskip

Let $Aut(\mathcal{T}^n)=\GLnZ$ denote the group of automorphisms of the group $\mathcal{T}^n=(\K^{*})^n$. As $\K^n$ contains $\mathcal{T}^n$ as an open subset, we get a natural injection  $T:\GLnZ\rightarrow \BirAn$ which image normalizes $\DKn$, the subgroup of $\GLn$ made up of \textit{diagonal automorphisms} $(x_1,...,x_n) \mapsto (\alpha_1 x_1,...,\alpha_n x_n)$, with $\alpha_i \in \K^{*}, i=1,...,n$. Identifying $\GLnZ$ with its image, the group of monomial birational maps with coefficients $1$, we will see $\GLnZ$ as a subgroup  of $\BirAn$. For a further description of the inclusion $\GLnZ \subset \BirAn$, see \cite{bib:GP}. \textit{(Note that there is another natural inclusion  from $\GLnZ$ to $\GLn \subset \BirAn$, but we won't use this one in this paper.)}

Among the affine transformations of $\An$, we distinguish those that we call \textit{almost-diagonal automorphisms}:  namely maps of the form $(x_1,...,x_n) \mapsto (x_1+1,\alpha_2 x_2,...,\alpha_n x_n)$. The $\alpha_i$ $\in \K^{*}$, $i=2,...,n$ will be called the \textit{eigenvalues} of the map. It is clear that $Aut(\mathcal{T}^{n-1})=\GLnmZ$ normalizes the set $\ADKn$ of almost-diagonal automorphisms.

We can now extend our diagram:
\begin{center}
$\begin{array}{cccccccc}
 &  \PGLnpl   & \subset &  & \BirPn \\
 &  \cup   & & & \| \\
  \GLn &\subset \AffKn \subset &\AutKn& \subset &\BirAn & \supset Aut(\mathcal{T}^n)&=&\GLnZ \\
  \cup &   \cup  & & & & \cup & & \cup \\
  \DKn& \ADKn    &  &  & &  Aut(\mathcal{T}^{n-1})&=&\GLnmZ.
\end{array}$
\end{center}
We will prove the following results:
\begin{theorem}[Conjugacy classes of affine automorphisms of $\An$]
\label{Thm:ResAn}\it
\begin{enumerate}
\item
{\it In the affine Cremona group  $\AutKn$ (Section \ref{Sec:ConjugationAutKn})}
\begin{itemize}
\item
The conjugacy classes of affine automorphisms that fix a point are given by their Jordan normal form, as in $\GLn$.
\item
Any affine automorphism that fix no point is conjugate to an almost-diagonal automorphism, unique up to a permutation of its eigenvalues.
\end{itemize}
\item
{\it In the Cremona group $\BirAn$ (Section \ref{Sec:ConjugationBirAn})}
\begin{itemize}
\item
Any affine automorphism of $\An$ is conjugate, either to a diagonal, or to an almost-diagonal automorphism of $\An$, exclusively.
\item
The conjugacy classes of diagonal and almost-diagonal automorphisms of $\An$ are respectively given by the orbits of the actions of $\GLnZ$ and $\GLnmZ$ by conjugation. These actions correspond respectively to the natural actions of $Aut(\mathcal{T}^n)$ and $Aut(\mathcal{T}^{n-1})$ on $\mathcal{T}^n$ and $\mathcal{T}^{n-1}$.
\item
In particular, if $n>1$, two affine automorphisms of the same finite order are conjugate in $\BirAn$.
\end{itemize}
\end{enumerate}
\end{theorem}

From theorem $\ref{Thm:ResAn}$ we can deduce the analogue result on $\Pn$, which is the following Theorem (Section \ref{Sec:ConjBirPn}). We define \textit{almost-diagonal automorphism of $\Pn$} to be a map of the form $(x_0:...:x_n) \mapsto (x_0:x_0+x_1:\alpha_2 x_2:...:\alpha_n x_n)$, with $\alpha_2,...,\alpha_n \in \K^{*}$.
\begin{theorem}[Conjugacy classes of automorphisms of $\Pn$ in the Cremona group $\BirPn$]
\it
\begin{itemize}
\item
Any automorphism of $\Pn$ is conjugate, either to a diagonal or to an almost-diagonal automorphism of $\Pn$, exclusively.
\item
The conjugacy classes of diagonal and almost-diagonal automorphisms of $\Pn$ are respectively given by the orbits of the actions of $\GLnZ$ and $\GLnmZ$ by conjugation. These actions correspond respectively to the natural actions of $Aut(\mathcal{T}^n)$ and $Aut(\mathcal{T}^{n-1})$ on $\mathcal{T}^n$ and $\mathcal{T}^{n-1}$.
\item
In particular, if $n>1$, two automorphisms of $\Pn$ of the same finite order are conjugate in $\BirAn$, (see \cite{bib:BBl}, Proposition 2.1.).
\end{itemize}
\end{theorem}
\section{Conjugacy classes of $\AffKn$}
Denote by $\K[X]=\K[X_1,...,X_n]$ the polynomial ring in the variables $X_1,...,X_n$ over $\K$ and by $\K(X)=\K(X_1,...,X_n)$ its field of fractions. Elements of $\BirAn$ can be written in the form $\varphi=(\varphi_1,...,\varphi_n)$, where each $\varphi_i$ belongs to $\K(X)$. Any birational map $\varphi\in \BirAn$ induces a
map \begin{center}$\varphi^{*}: F \mapsto F \circ \varphi$, ($F \in \K(X)$)\end{center}
which is a $\K$-automorphism of the field $\K(X)$. Conversely, any $\K$-automorphism of $\K(X)$ is of this form. So, $\BirAn$ is anti-isomorphic to the group of $\K$-automorphisms of $\K(X)$ and
its subgroups $\AutKn$, $\AffKn$ and $\GLn$ corresponds respectively to the groups of $\K$-automorphisms of $\K[X]$, $\K[X]_{\leq 1}$ and $\K[X]_{1}$. Here, $\K[X]_{\leq 1}$ and $\K[X]_{1}$ denote respectively the sets of polynomials of degree~$\leq 1$ and egal to $1$.

\label{Sec:ConjugationAffAn}
The study of the conjugacy classes of $\AffKn$ is elementary and well-known. Let $\alpha,\beta \in \AffKn$; 
let us recall that the first dichotomy consists in separating the cases according to whether $\alpha$ and $
\beta$ fix a point or not, since an affine automorphism that fixes a point cannot be conjugate to one with no fixed point.

If both $\alpha$ and $\beta$ fix a point, they are respectively conjugate to linear automorphisms $\alpha'$ and $\beta'$ of $\K^n$ (elements of $\GLn$), and are then conjugate if and only if these have the same Jordan normal form. We will say that these Jordan normal forms are also the \textit{Jordan normal forms} of $\alpha$ and $\beta$. That doesn't depend on the choice of $\alpha'$ and $\beta'$.

We will extend this idea to the case of affine automorphisms with no fixed point. Suppose that $\alpha$ has no fixed point. We consider a basis $(1,P_1,...,P_n)$ of $\K[X]_{\leq 1}$ such that the matrix of ${\alpha^{*}}_{|_{\K[X]_{\leq 1}}}$ has the Jordan normal form 
$\left(\begin{array}{ccc} J(\lambda_1,k_1)& & \\  & \ddots & \\  & &J(\lambda_r,k_r)\end{array}\right)
$, where  $J(\mu,k) = \left(\begin{array}{cccc} \mu & 1 & & \\ & \ddots & \ddots & \\ & & \mu & 1  \\ & &  &\mu \end{array}\right)\in GL(k,\K)$ is a Jordan block of size $k$. Observe that $\lambda_1=1$ and $k_1>1$, as $\alpha$ fixes no point.

Then, $\pi:(x_1,...,x_n) \mapsto (P_1(x_1,...,x_n),...,P_n(x_1,...,x_n))$ is an affine automorphism of $\An$ such that 
\begin{equation}\label{JordanAffine}
\pi \alpha \pi^{-1}:\left(\begin{array}{c}x_1 \\x_2 \\  \vdots \\ x_n\end{array}\right) \mapsto
 \leftexp{t}{\left(\begin{array}{ccc} J(\lambda_1,k_1-1)& & \\  & \ddots & \\  & &J(\lambda_r,k_r)\end{array}\right)}
\left(\begin{array}{c}x_1 \\ x_2 \\ \vdots \\ x_n\end{array}\right)
+\left(\begin{array}{c}1 \\ 0 \\ \vdots \\ 0\end{array}\right);
\end{equation}
by analogy, we will also say that $(\ref{JordanAffine})$ is the \textit{Jordan normal form of $\alpha$}.

By observing that an automorphism given by formula $(\ref{JordanAffine})$ is of infinite order (since $char(\K)=0$), we see that any affine automorphism of finite order has a fixed point, which is true in general for any automorphism of finite order (see \cite{bib:KS}).

Looking at the linear action of $\alpha^{*}$ on $\K[X]_{\leq 1}$, we observe that the pairs $(\lambda_i,k_i)$ characterize the conjugacy class of $\alpha$.
We have thus showed the following proposition:

\begin{proposition}[Conjugacy classes in $\AffKn$]
\label{Prp:ConjugacyAffAn}
{\it Two affine automorphisms of $\An$ are conjugate if and only if they have the same Jordan normal form.}
\proofend
\end{proposition}
\section{Conjugation in the group $\AutKn$}
\label{Sec:ConjugationAutKn}
It is clear that the conjugation in $\AutKn$ respects the dichotomy of the existence or not of a fixed point.

When there exists a fixed point, passing from the group $\AffKn$ to $\AutKn$ does not bring anything new:

\bigskip

\begin{proposition}
\label{Prp:GLnAutKn}
{\it Two affine automorphisms of $\An$ that fix a point are conjugate in $\AutKn$ if and only if they are already conjugate in $\AffKn$.}\end{proposition}
\begin{proof}
Let $\alpha$, $\beta$ be two affine automorphisms which have fixed points. Changing $\alpha$ and $\beta$ within their $\AffKn$-conjugacy classes, we can suppose that $\alpha$ and $\beta$ belong to $\GLn$. Let $\pi \in \AutKn$ be such that $\pi \alpha=\beta \pi$ and let $\rho \in \GLn$ denote the tangent map of $\pi$ at the origin; then $\rho\alpha=\beta\rho$.
\proofend\end{proof}

On the other hand, if there are no fixed points, by means of elements of $\AutKn$, it is possible to modify the size of the Jordan blocks as the following shows:

\bigskip

\begin{example}
\label{exa1}
{\it The affine automorphisms
\begin{center}
$\alpha:\left(\begin{array}{c}x_1  \\  x_2\end{array}\right) \mapsto 
\left(\begin{array}{cc}1 & 0 \\  1 & 1\end{array}\right)\left(\begin{array}{c}x_1  \\  x_2\end{array}\right)+\left(\begin{array}{c}1  \\  0\end{array}\right)$
\\
$\beta:\left(\begin{array}{c}x_1  \\  x_2\end{array}\right) \mapsto 
\left(\begin{array}{cc}1 & 0 \\  0 & 1\end{array}\right)\left(\begin{array}{c}x_1  \\  x_2\end{array}\right)+\left(\begin{array}{c}1  \\  0\end{array}\right)$.\end{center}

are conjugate by the automorphism  $\pi:\left(\begin{array}{c}x_1  \\  x_2\end{array}\right) \mapsto \left(\begin{array}{c}x_1  \\  x_2-\frac{x_1(x_1-1)}{2}\end{array}\right)$.

Geometrically, the affine automorphisms $\alpha$ and $\beta$ on $\K^2$ respectively leave invariant the conics $C_t$ given by $x_2+\frac{x_1-{x_1}^2}{2}=t$ and the lines $L_t$ given by $x_2=t$ where $t \in \K$; the action of these automorphisms on this curves is just a translation by $1$ on the axis $x_1$. The automorphism $\pi$ sends every conic $C_t$ on the line $L_t$ without changing the coordinate $x_1$, and so sends the orbits of $\alpha$ on those of $\beta$.}
\end{example}

\begin{example}
{\it The affine automorphisms
\begin{center}
$\begin{array}{cccc}
\alpha:\left(\begin{array}{c}x_1 \\ x_2 \\ \vdots  \\ x_n\end{array}\right) &\mapsto &
\left(\begin{array}{cccc} 1 & & & \\ 1& 1 & & \\ & \ddots & \ddots &   \\ & & 1 &1 \end{array}\right)&\left(\begin{array}{c}x_1 \\ x_2 \\ \vdots  \\ x_n\end{array}\right)+\left(\begin{array}{c}1 \\ 0 \\ \vdots  \\ 0\end{array}\right)
\\
\beta:\left(\begin{array}{c}x_1 \\ x_2 \\ \vdots  \\ x_n\end{array}\right) &\mapsto &
\left(\begin{array}{cccc} 1 & & & \\ & 1 & & \\ &  & \ddots &   \\ & &  &1 \end{array}\right)&\left(\begin{array}{c}x_1 \\ x_2 \\ \vdots  \\ x_n\end{array}\right)+\left(\begin{array}{c}1 \\ 0 \\ \vdots  \\ 0\end{array}\right)\end{array}$\end{center}

are conjugate by the automorphism  
\begin{center}$\pi:\left(\begin{array}{c}x_1 \\ x_2 \\ x_3\\\vdots  \\ x_n\end{array}\right) \mapsto 
\left(\begin{array}{l}x_1 \\ x_2+P_2(x_1) \\ x_3+P_3(x_1,x_2) \\ \vdots  \\ x_n+P_n(x_1,...,x_{n-1})\end{array}\right),$\end{center}
where $P_m \in \K[X_1,...,X_{m-1}]$ is defined by the formula $$P_m=\sum_{k=1}^{m-2} \left(\begin{array}{c} x_1+k-1 \\ k\end{array}\right) (-1)^k x_{m-k}+(-1)^{m-1} (m-1) \left(\begin{array}{c} x_1+m-2 \\ m\end{array}\right),$$
where we denote by $\left(\begin{array}{c} Q \\ r\end{array}\right)$ the polynomial $\frac{1}{r!}Q(Q-1)(Q-2)\cdots (Q-(r-1))$, for $r \in \mathbb{N}, Q \in \K[X_1,...,X_n]$.

We give a more precise idea, we explicit the polynomials $P_2,...,P_5$:

$\begin{array}{lll}
P_2&=&-\frac{1}{2}x_1(x_1-1)  \mbox{(see Example \ref{exa1})}\\
P_3&=&-x_1x_2+\frac{1}{3}(x_1-1)x_1(x_1+1)\\
P_4&=&-x_1x_3+\frac{1}{2}x_1(x_1+1)x_2-\frac{1}{8}(x_1-1)x_1(x_1+1)(x_1+2)\\
P_5&=&-x_1x_4+\frac{1}{2}x_1(x_1+1)x_3-\frac{1}{6}x_1(x_1+1)(x_1+2)x_2\\
& &\hspace{4 cm} -\frac{1}{30}(x_1-1)x_1(x_1+1)(x_1+2)(x_1+3)\end{array}$

Geometrically, the affine automorphisms $\alpha$ and $\beta$ on $\K^2$ respectively leave invariant the curves given by $x_2+P_2(x_1)=\tau_2,x_3+P_3(x_1,x_2)=\tau_2$,...,$x_n+P_{n}(x_1,...,x_{n-1})=\tau_n$ and the lines given by $x_i=\tau_i, i=2,..,n$ where $(\tau_2,...,\tau_n) \in \K^{n-1}$; the action of these automorphisms on this curves is just a translation by $1$ on the axis $x_1$. The automorphism $\pi$ sends the curves on the lines without changing the coordinate $x_1$, and so sends the orbits of $\alpha$ on those of $\beta$.}
\end{example}

\bigskip

In fact, we can generalize this example to reduce a lot the size of the Jordan blocks:

\begin{proposition}
\label{Prp:ConjugAlmostDiagAutAn}
{\it An affine automorphism of $\An$ with no fixed point  is conjugate, in $\AutKn$, to an almost-diagonal automorphism, unique up to a permutation of its eigenvalues.

More precisely, the automorphism 

\begin{center}$\alpha:\left(\begin{array}{c}x_1 \\x_2 \\  \vdots \\ x_n\end{array}\right) \mapsto
 \leftexp{t}{\left(\begin{array}{cccc} J(\lambda_1,k_1-1)& & \\ & J(\lambda_2,k_2) \\ &  & \ddots & \\  & & &J(\lambda_r,k_r)\end{array}\right)}
\left(\begin{array}{c}x_1 \\ x_2 \\ \vdots \\ x_n\end{array}\right)
+\left(\begin{array}{c}1 \\ 0 \\ \vdots \\ 0\end{array}\right)$, with $\lambda_1=1$,
\end{center}
on Jordan normal form, is conjugate, in $\AutKn$, to the almost-diagonal automorphism
\begin{center}
$\alpha_{D}:\left(\begin{array}{c}x_1 \\x_2 \\  \vdots \\ x_n\end{array}\right) \mapsto
{\left(\begin{array}{cccc} \lambda_1\cdot Id_{k_1-1}& & \\ & \lambda_2\cdot Id_{k_2}& & \\&  & \ddots & \\ & & &\lambda_r \cdot Id_{k_r}\end{array}\right)}
\left(\begin{array}{c}x_1 \\ x_2 \\ \vdots \\ x_n\end{array}\right)
+\left(\begin{array}{c}1 \\ 0 \\ \vdots \\ 0\end{array}\right)$,\end{center}where $\lambda \cdot Id_k=\left(\begin{array}{ccc}\lambda \\ & \ddots \\ & & \lambda\end{array}\right) \in GL(k,\K)$ is the diagonal part of $J(\lambda,k)$.}
\end{proposition}
\begin{proof}\textit{   }
\begin{enumerate}
\item
{\it The conjugation} \\
For an automorphism $\pi:(x_1,...,x_n) \mapsto (P_1(x_1,...,x_n),...,P_n(x_1,...,x_n))$, we have $\pi \circ \alpha=\alpha_D \circ \pi$ if and only if
\begin{equation}
\alpha^{*}(P_j)=\left\{
\begin{array}{llll}
P_1+1 &\mbox{if }j=1  \\ 
\mu(j) P_j& \mbox{if  } j>1
\label{Condition:Conjugate}
\end{array}\right.\end{equation}
where $\mu(j)$ denotes the $j$-th eigenvalue of $\alpha_D$.

If $j$ is the first indice of the $i$-th block, set $P_j=X_j$.
If not, there exists $Q_j \in \K[X_1,...,X_{j-1}]$ such that $\alpha^{*}(Q_j)=\mu(j)Q_j-X_{j-1}$ (see Lemma $\ref{Lem:techsurjective}$ below); then let $P_j=X_j+Q_j$ so that $\alpha^{*}(P_j)=\alpha^{*}(X_j)+\alpha^{*}(Q_j)=(\mu(j)X_j+X_{j-1})+(\mu(j)Q_j-X_{j-1})=\mu(j)(X_j+Q_j)=\mu(j) P_j$. It is clear that the map defined by these $P_j$ is an automorphism of $\An$ that conjugates $\alpha$ and $\alpha_D$.
\item
{\it The unicity} \\
Let us suppose that the two almost-diagonal automorphisms
\begin{center}$\begin{array}{cccccc}
\theta_{\mu}:&(x_1,...,x_n)& \mapsto& (x_1,\mu_2 x_2,\mu_3 x_3,...,\mu_n x_n)&+&(1,0,...,0)\\
\theta_{\nu}:&(x_1,...,x_n) &\mapsto &(x_1,\nu_2 x_2,\nu_3 x_3,...,\nu_n x_n)&+&(1,0,...,0)\end{array}$\end{center}
are conjugate in $\AutKn$: there exists $\pi \in \AutKn$ such that $\theta_{\mu}\circ \pi=\pi \circ \theta_{\nu}$.
By derivation, one sees that 
\begin{equation}\label{TwoAlmostUnicity}
\left(\begin{array}{cccc}
1 \\ & \mu_2 \\ & & \ddots \\ & & & \mu_n\end{array}\right) \circ M(x) = M(\theta_{\mu}(x)) \circ \left(\begin{array}{cccc}
1 \\ & \nu_2 \\ & & \ddots \\ & & & \nu_n\end{array}\right)\end{equation}
where $M(x)=(m_{ij}(x))_{i,j=1}^n$ denotes the Jacobian matrix of $\pi$ at $x=(x_1,...,x_n)$.

So we have $\mu_i m_{ij}(x) \nu_j^{-1}=m_{ij}(x_1+1,\nu_2 x_2,...,\nu_n x_n)$, where $\mu_1=\nu_1=1$. Hence the $m_{ij}$ are eigenvectors of $\theta_{\nu}^{*}$ and belong to $\K[X_2,...,X_n]$, by Lemma $\ref{Lem:NoX1OnEigenVect}$ below. In particular, we see that $M(0,...,0)=M(1,0,...,0)=M(\theta_{\nu}(0,...,0))$. Since $\pi$ is an automorphism of $\K^n$, we have $det(M(x))=det(M(0,...,0))\in \K^{*}$, so equality ($\ref{TwoAlmostUnicity}$) evaluated at $x=(0,...,0)$ shows that the diagonal matrices 
\begin{center}
$\left(\begin{array}{cccc}
1 \\ & \mu_2 \\ & & \ddots \\ & & & \mu_n\end{array}\right)$ and $\left(\begin{array}{cccc}
1 \\ & \nu_2 \\ & & \ddots \\ & & & \nu_n\end{array}\right)$ 
\end{center}
are conjugate in $GL(n,\K)$, i.e. the $\mu_i$ and the $\nu_i$ are equal up to permutation.
\proofend
\end{enumerate}
\end{proof}

\bigskip

\begin{remark}
\begin{enumerate}\item
In fact, Proposition \ref{Prp:ConjugAlmostDiagAutAn} shows also that $\alpha$ and $\alpha_D$ are conjugate in the famous \textit{Jonqui\`ere subgroup} of $\AutKn$ given by $\{\varphi=(\varphi_1,...,\varphi_n) \in \AutKn \ | \ \varphi_i \in \K[X_1,...,X_i] \mbox{ for } i=1,...,n\}$.
\item
The characteristic $0$ is very important here, because Lemma $\ref{Lem:techsurjective}$ (and then Proposition \ref{Prp:ConjugAlmostDiagAutAn}) is false in characteristic $>0$:
let  $\mathbb{L}$ be any field of characteristic $p$, and $\alpha_1,\alpha_2,...,\alpha_p$ be  the affine automorphisms \begin{center}$\alpha_k:(x_1,x_2,...,x_n) \mapsto (x_1+1,x_2+x_1,...,x_k+x_{k-1},x_{k+1},x_{k+2},...,x_n)$.\end{center} Then $\alpha_1,\alpha_2,...,\alpha_{p-1}$ are all of order $p$ and conjugate in $Bir(\mathbb{L}^n)$ but not $\alpha_p$ which is of order $p^2$.
\end{enumerate}
\end{remark}

\bigskip

\begin{lemma}
\label{Lem:techsurjective}
{\it If $j$ is not the first indice of a block, there exists $Q_j \in \K[X_1,...,X_{j-1}]$ such that $\alpha^{*}(Q_j)=\mu(j)Q_j - X_{j-1}$.}
\end{lemma}
\begin{proof}
Let us recall that :
\begin{center}$\alpha^{*}(X_k)=\left\{
\begin{array}{ll}X_1+1 &\mbox{if $k$ is the first indice}\\
 & \ \mbox{of the first block $(k=1)$}\\
\lambda(k)X_k & \mbox{if $k$ is the first indice of another block}\\
\lambda(k)X_k+X_{k-1} & \mbox{if $k$ is not the first indice of a block}\end{array}\right.$\end{center}

We will prove the stronger assertion (that implies the Lemma, using $k=j-1$, as $\mu(j-1)=\mu(j)$):
\begin{equation}\label{assertionX1t}
\begin{array}{c}\mbox{ \it For any integers }1\leq k \leq n\mbox{ \it and }t\geq 0, \mbox{\it the monomial } X_1^tX_k  \mbox{\it \ is } \\ \mbox{\it in the image of the linear map of vector spaces } \\ (\alpha^{*}-\mu(k)id):\K[X_1,...,X_k] \rightarrow \K[X_1,...,X_k].
\end{array}\end{equation}

\begin{itemize}
\item
Since $(\alpha^{*}-id)_{|_{\K[X_1]}}$ is surjective, the assertion $(\ref{assertionX1t})$ is true for $k=1$ and $t\geq 0$.
\item
In the same way, if $k>1$ is the first indice of a block, since $\alpha^{*}(X_k)=\mu(k)X_k$, the linear map $(\alpha^{*}-\mu(k) id)_{|_{\K[X_1]\cdot X_k}}$ is surjective and the assertion $(\ref{assertionX1t})$ is correct for $k$ and $t\geq 0$.
\item
Lastly, if $k$ is not the first indice of a block, we have 
\begin{center}$\begin{array}{ccl}
\alpha^{*}(X_1^{t+1}X_k)&=&(X_1+1)^{t+1}(\mu(k)X_k+X_{k-1})\\
&=&\mu(k)X_1^{t+1}X_k+(t+1)\mu(k)X_1^t X_k+\sum_{l,s}a_{l,s} X_1^s X_k,\end{array}$\end{center}
where $a_{l,s} \in \K$, and all the $(l,s)$ are strictly smaller than $(k,t)$, for the lexicographical order.
\end{itemize}
The assertion $(\ref{assertionX1t})$ is then right by induction on the indices $(k,t)$.
\proofend
\end{proof}

\bigskip

\begin{lemma}
\label{Lem:NoX1OnEigenVect}
{\it For any almost-diagonal automorphism $\theta_{\lambda}:(x_1,x_2,...,x_n) \mapsto(x_1+1,\lambda_2 x_2,...,\lambda_{n} x_{n})$  of $\An$, the eigenvectors of $\theta_{\lambda}^{*}$ belong to $\K[X_2,...,X_n]$.}
\end{lemma}\begin{proof}
We first observe that the map $\theta_{\lambda}^{*}$ leaves invariant the decomposition $$\K[X]= \bigoplus_{(a_2,...,a_n) \in \mathbb{N}^{n-1}} (\prod_{i\geq 2}^n {X_i}^{a_i} \cdot \K[X_1])$$ and that the $\Pi_{i\geq 2} {X_i}^{a_i}$ are eigenvectors. Since the only eigenspace of $\theta_{\lambda}^{*}$ in $\K[X_1]$ is $\K\cdot 1$, any eigenvector is in $\K[X_2,...,X_n]$.
\proofend
\end{proof}
\section{Conjugation in the group $\BirAn$}
\label{Sec:ConjugationBirAn}
\subsection{Diagonalizable and birationaly almost-diagonal affine automorphisms}
The conjugation of elements of $\AffKn$ in $\BirAn$ changes the dichotomy on the existence of a fixed point, as the following example shows:

\bigskip

\begin{example}
\label{ExampleJordAlmDiag}
\textit{The two affine automorphisms
\begin{center}
$\begin{array}{lll}
\alpha:\left(\begin{array}{c}x_1  \\  x_2\end{array}\right) &  \mapsto &
\left(\begin{array}{cc}1 & 0 \\  1 & 1\end{array}\right)\left(\begin{array}{c}x_1  \\  x_2\end{array}\right) \\
\beta:\left(\begin{array}{c}x_1  \\  x_2\end{array}\right)& \mapsto &
\left(\begin{array}{cc}1 & 0 \\  0 & 1\end{array}\right)\left(\begin{array}{c}x_1  \\  x_2\end{array}\right)+\left(\begin{array}{c}1  \\  0\end{array}\right)
\end{array}$\end{center}
are conjugate by the birational map  \begin{center}$\varphi:\left(\begin{array}{c}x_1  \\  x_2\end{array}\right) \mapsto \left(\begin{array}{c}x_2x_1^{-1}  \\  x_1^{-1}\end{array}\right).$\end{center}
In fact, choosing $\K^2$ as the open subset ${x_0\not=0}$ of $\mathbb{P}^2$, the two affine automorphisms become
\begin{center}$\begin{array}{ccc}\tilde{\alpha}:(x_0:x_1:x_2) & \mapsto & (x_0:x_1:x_2+x_1)\\
\tilde{\beta}:(x_0:x_1:x_2) & \mapsto & (x_0:x_1+x_0:x_2)\end{array}$\end{center}
and the birational map corresponds only to the permutation of coordinates \begin{center}
$\tilde{\varphi}:(x_0:x_1:x_2) \mapsto (x_1:x_2:x_0)$.\end{center}}
\end{example}

\bigskip

Let us recall that $\DKn$ denotes the subgroup of $\GLn$ made up of diagonal automorphisms and $\ADKn$ the subset of $\AffKn$ made up of almost-diagonal automorphisms. We will say that $\alpha \in \AffKn$ is \textit{diagonalizable} (respectively \textit{birationaly almost-diagonal}) if there exists $\pi \in \AffKn$ such that $\pi\alpha\pi^{-1} \in \DKn$ (respectively if there exists $\pi \in \BirAn$ such that $\pi\alpha\pi^{-1} \in \ADKn$).

These two notions provide the dichotomy for conjugation in $\BirAn$, as the following result shows:

\bigskip

\begin{proposition}
{\it An affine  automorphism of $\An$ is either diagonalizable or birationaly almost-diagonal, exclusively.}
\end{proposition}
\begin{proof}\textit{   }
\begin{enumerate}
\item
{\it Existence of the conjugation} \\
Let $\alpha \in \AffKn$. If $\alpha$ fixes no point, then $\alpha$ is birationaly almost-diagonal, by Proposition ($\ref{Prp:ConjugAlmostDiagAutAn}$). Else, by a conjugation in $\AffKn$, we may suppose that $\alpha \in \GLn$; if $\alpha$ is not diagonalizable, we verify, as in the example $(\ref{ExampleJordAlmDiag})$ above, that $\alpha$ is birationaly almost-diagonal.

Explicitly, let 
$\alpha$ be of the form $x \mapsto 
\leftexp{t}{\left(\begin{array}{c|c} J(\mu,m_1)& \\ \hline   & \ddots
\end{array}\right)}x$, with $m_1>1$. Then, the birational map \begin{center}$\pi:(x_1,...,x_n) \dasharrow (\frac{1}{x_1},\frac{x_2}{x_1},...,\frac{x_{m_1}}{x_1},x_{m_1+1},x_{m_1+2},...,x_n)$\end{center} conjugates $\alpha$ to the affine automorphism 
\begin{center}
$x \mapsto 
{\left(\begin{array}{ccccc|c} \mu^{-1} & & & &  \\ & 1& & & &  \\ & \mu^{-1} & 1 & & & \\ & & \ddots & \ddots & &  \\ & & &  \mu^{-1} & 1 & \\ \hline   & & & & &\ddots
\end{array}\right)}x+ \left(\begin{array}{c}0 \\  \mu^{-1} \\0 \\ \vdots \\ \vdots \\ 0 \end{array}\right)$,\end{center}
that fixes no point and so that is  birationally almost-diagonal. \\
\item
{\it The diagonal and almost-diagonal automorphisms of $\An$ are in distinct conjugacy classes of~$\BirAn$}\\
Suppose that $\theta_{\lambda}$ is an almost-diagonal automorphism of $\An$ conjugate in $\BirAn$ to a diagonal automorphism $\rho_{\mu}$ by a birational map $\pi$. We write $\rho_{\mu}, \theta_{\lambda}, \pi$ in their explicit form:
\begin{center}$\begin{array}{ccccccccccccc}
\theta_{\lambda}&:&(x_1,x_2,...,x_n) &\mapsto&(&x_1+1&,&\lambda_2 x_2&,&...&,&\lambda_{n} x_{n}&)\\
\rho_{\mu}&:& (x_1,x_2,...,x_n) &\mapsto&  (&\mu_1 x_1&,&\mu_2 x_2&,&...&,&\mu_{n} x_{n}&)\\
\pi&:& (x_1,x_2,...,x_n)& \mapsto  &(&\frac{P_1(x_1,...,x_n)}{Q_1(x_1,...,x_n)}&,&\frac{P_2(x_1,...,x_n)}{Q_2(x_1,...,x_n)}&,&...&,&\frac{P_n(x_1,...,x_n)}{Q_n(x_1,...,x_n)}&)
\end{array}$
\end{center}
with $\mu_i,\lambda_i \in \K^{*}$ and $P_i,Q_i \in \K[X]$ without common divisors. Since $\theta_{\lambda}^{*}$ gives an automorphism of $\K[X]$ the condition $\pi \theta_{\lambda}=\rho_{\mu} \pi$ implies that $\frac{\theta_{\lambda}^{*}(P_i)}{\theta_{\lambda}^{*}(Q_i)}=\mu_i \frac{P_i}{Q_i}$, for any $i=1,...,n$. Then all the $P_i$ and $Q_i$ must be eigenvectors of $\theta_{\lambda}^{*}$, viewed as a $\K$-linear map. But such eigenvectors belongs to $\K[X_2,...,X_n]$ (Lemma $\ref{Lem:NoX1OnEigenVect}$) and so the map $\pi$ cannot be birational.
\proofend
\end{enumerate}
\end{proof}

\bigskip

\begin{remark}
\textit{ The fact that a non-diagonalizable automorphism of $\K^n$ is conjugate, in $\BirAn$, to an affine automorphism with no fixed point can also be viewed as follows:}

Let us consider the automorphism
\begin{center}$\alpha:\left(\begin{array}{c}x_1  \\ \vdots \\ x_n \end{array}\right) \mapsto 
{\left(\begin{array}{cccc|c} \mu & & &  \\ 1& \mu&  & &  \\ & \ddots & \ddots & &   \\ & & 1 & \mu & \\ \hline   & & & & \ddots
\end{array}\right)}\left(\begin{array}{c}x_1  \\ \vdots \\ x_n \end{array}\right)$\end{center} of $\K^n$, and his extension to an automorphism
\begin{center}$\tilde{\alpha}:\left(\begin{array}{c}x_0 \\ x_1  \\ \vdots \\ x_n \end{array}\right)\mapsto 
{\left(\begin{array}{c|cccc|c}1 & & & & & \\ \hline &   \mu & & & &  \\ &1& \mu&  & &  \\ & & \ddots & \ddots & &   \\ & & & 1 & \mu & \\ \hline   & & & & & \ddots
\end{array}\right)}\left(\begin{array}{c}x_0 \\ x_1  \\ \vdots \\ x_n \end{array}\right)$\end{center} of $\mathbb{P}^n(\K)$. We observe that $\tilde{\alpha}$ leaves invariant the open set $(x_1 \not=0)$ of $\mathbb{P}^n(\K)$ and it induces in it the affine automorphism
\begin{center}$\beta:\left(\begin{array}{c}x_0  \\ x_2 \\ \vdots \\ x_n \end{array}\right) \mapsto 
{\left(\begin{array}{ccccc|c} \mu^{-1} & & & &  \\ & 1&  & & &  \\ &  \mu^{-1} & 1 & &  \\ & & \ddots & \ddots  & &  \\ & & &\mu^{-1} & 1 &   \\ \hline   & & & & & \ddots 
\end{array}\right)}\left(\begin{array}{c}x_0 \\ x_2 \\ x_3 \\ \vdots \\ x_n \end{array}\right)+\left(\begin{array}{c}0 \\ 1 \\ 0 \\ \vdots \\ 0 \end{array}\right)$\end{center} 
of $\K^n$, that fixes no point. We go from $\alpha$ to $\beta$ by exchange of $x_0$ and $x_1$.
\end{remark}

\subsection{Actions of $\GLnZ$ and $\GLnmZ$}
\label{subsection:Actions}
Let us summarize the situation;  we have
\begin{center}
$\AffKn= \AffKn \cap (\BirAn \bullet \DKn \uplus \BirAn \bullet \ADKn)$
\end{center}
where $\BirAn \bullet B= \{ \pi \alpha \pi^{-1} \ | \ \alpha \in B, \pi \in \BirAn\}$.

\bigskip

We continue our study by describing the trace on $\DKn$ (respectively on $\ADKn$) of the conjugacy class of an element $\rho_{\mu}$ (respectively $\theta_{\nu}$) of $\DKn$ (respectively $\ADKn$). Equivalently, we are looking for:
\begin{center}
$\begin{array}{ccc}
(\BirAn \bullet \rho_{\mu}) &\cap& \DKn \\
(\BirAn \bullet \theta_{\nu}) &\cap& \ADKn\end{array}$ \end{center}

{\bf Action of $\GLnZ$ on $\DKn$}\\
For this purpose, let us recall that $\mathcal{T}^n=(\K^{*})^n$ and consider the isomorphism
\begin{center}$\begin{array}{cccc}
\rho: & \mathcal{T}^n & \rightarrow & \DKn \subset \BirAn \\
              & (\mu_1,....,\mu_n) & \mapsto & [(x_1,...,x_n) \mapsto (\mu_1 x_1,...,\mu_n x_n)]\end{array}$\end{center}
that is $\GLnZ$-equivariant for actions that we discribe now:

The action on $\mathcal{T}^n$ is the natural action of  of $\GLnZ=Aut(\mathcal{T}^n)$: a matrix $A=(a_{ij})_{i,j=1}^n \in \GLnZ$ maps an element $(\mu_1,....,\mu_n) \in (\K^{*})^n$ on \begin{center}$A \bullet (\mu_1,...,\mu_n) = (\mu_1^{a_{11}}\mu_2^{a_{12}}\cdots\mu_n^{a_{1n}},\mu_1^{a_{21}}\mu_2^{a_{22}}\cdots\mu_n^{a_{2n}},...,\mu_1^{a_{n1}}\mu_2^{a_{n2}}\cdots\mu_n^{a_{nn}})$.\end{center}
On the other side, we have the  injective homomorphism
\begin{center}$\begin{array}{cccc}
T: & \GLnZ & \rightarrow & \BirAn \\
 & A=(a_{ij})_{i,j=1}^n & \mapsto & [(x_1,...,x_n) \dasharrow (x_1^{a_{11}}\cdots x_n^{a_{1n}},...,x_1^{a_{n1}}\cdots x_n^{a_{nn}})]\end{array}$\end{center}
whose image normalizes $\DKn$. This gives then an action of $\GLnZ$ on $\DKn$ by conjugation.

We get the formula:
\begin{equation}
   \begin{array}{|c|}
\hline
T(A) \circ \rho(\mu) \circ {T(A)}^{-1}=\rho(A \bullet\mu)\\
\hline
   \end{array} \ \ \mu \in (\K^{*})^n, A \in \GLnZ 
   \label{Condition:ActionGLn}\end{equation}
   that proves that $\rho$ is equivariant and that
\begin{center}   $(\BirAn \bullet \rho(\mu) )\cap \DKn \supset \rho(\GLnZ \bullet \mu)$.\end{center}
We will see further that this is an equality, i.e. that two diagonal automorphisms of $\An$ that are conjugate in $\BirAn$ are conjugate by an element of $\GLnZ$ (see Proposition \ref{Prp:GroupActionDiag}).

\bigskip

{\bf Action of $\GLnmZ$ on $\ADKn$}\\
We do the same for $\ADKn$: we have the injective homomorphism
\begin{center}$\begin{array}{cccc}
S: & \GLnmZ & \rightarrow & \BirAn \\
 & A=(a_{ij})_{i,j=2}^n & \mapsto & [(x_1,...,x_n) \dasharrow (x_1,x_2^{a_{21}}\cdots x_n^{a_{2n}},...,x_2^{a_{n2}}\cdots x_n^{a_{nn}})]\end{array}$\end{center}
whose image normalizes $\ADKn$ and who gives then an action of $\GLnmZ$ on $\ADKn$ by conjugation.

With the action of $\GLnmZ$ on $(\K^{*})^{n-1}$, the bijection
\begin{center}$\begin{array}{cccc}
\theta: & (\K^{*})^{n-1} & \rightarrow & \ADKn \subset \BirAn \\
              & (\nu_2,....,\nu_n) & \mapsto & [(x_1,...,x_n) \mapsto (x_1+1,\nu_2 x_2...,\nu_n x_n)]\end{array}$\end{center}
is $\GLnmZ$-equivariant, by the formula
\begin{equation}
   \begin{array}{|c|}
\hline
S(B) \circ \theta(\nu) \circ {S(B)}^{-1}=\theta(B \bullet\mu)\\
\hline
   \end{array} \ \ \nu \in (\K^{*})^{n-1}, B \in \GLnmZ 
   \label{Condition:ActionGLnm}\end{equation}
   and we see too that
   \begin{center}   $(\BirAn \bullet \theta(\nu) )\cap \ADKn \supset \theta(\GLnmZ \bullet \nu)$.\end{center}
We will see further that this is an equality, i.e. that two almost-diagonal automorphisms of $\An$ that are conjugate in $\BirAn$ are conjugate by an element of $\GLnmZ$ (see Proposition \ref{Prp:Proof:Thm: GroupActionSemiDiag}).

\subsection{Elements of finite order}
Let us prove the equality between $(\BirAn \bullet \rho(\mu) )\cap \DKn$ and $\rho(\GLnZ \bullet \mu)$, for elements $\mu \in (\K^{*})^n$ of finite order. For $n>1$, we do this by showing that $\rho(\GLnZ \bullet \mu)$ contains all diagonal automorphisms of the same order as $\mu$. This also proves that two affine automorphisms of the same finite order are conjugate in $\BirAn$.

If $n=1$ the group $Bir(\K)$ is equal to $PGL(2,\K)$ and a simple calculation shows that two diagonal automorphisms $x \mapsto \alpha x$ and $x \mapsto \beta x$ are conjugate if and only if $\alpha=\beta^{\pm 1}$. So the equality is true in dimension one too.
\begin{proposition}[Diagonal automorphisms of finite order]
\label{Prp:FiniteOrder}
{\it 
If $n>1$, two diagonal automorphisms of $\An$ of the same (finite) order are conjugate by an element of $\GLnZ \subset \BirAn$.}
\end{proposition}
\begin{proof}
Let $m$ be a positive integer and $\xi$ a primitive $m$-th root of the unity. A diagonal automorphism of $\An$ of order $m$ is given by $\rho(\alpha):(x_1,...,x_n) \mapsto (x_1 \xi^{t_1},...,x_n \xi^{t_n})$, where $\alpha=(\xi^{t_1},...,\xi^{t_n})$, and the greatest common divisor of the exponents $g=gcd(t_1,...,t_n)$ is prime to $m$.

The action of $\GLnZ$ on the diagonal automorphisms corresponds here to the usual action of $\GLnZ$ on the exponents $(t_1,...,t_n) \in \mathbb{Z}^n$.

With elementary matrices, that add a multiple of a coordinate to another one, we can map the vector $t=(t_1,...,t_n)$ to the vector $(g,0,...,0)$, so the diagonal automorphism $\rho(\alpha)$ is in the same orbit as the automorphism $(x_1,...,x_n)\mapsto (x_1 \xi^{g},x_2,...,x_n)$.
Because $g$ is prime to $m$, there exist $p,q \in \mathbb{Z}$ such that $pm+gq=1$. Since
 \begin{center}$ 
  \left(
\begin{array}{ccccc}
  q & -p&  & &   \\
  m & g &  & & \\
      &   & 1 &  & \\
      &   & & \ddots  & \\
      &   &  &  &1 \\
\end{array}
\right)
   \left(
\begin{array}{c}
 g \\
0\\
  \vdots\\
  \vdots\\
  0
\end{array}
\right) =
   \left(
\begin{array}{c}
 gq  \\
gm \\
 0\\
  \vdots\\
  0
\end{array}
\right) $
  \end{center}
we see that $\rho(\alpha)$ is conjugate, by an element of $\GLnZ$, to the diagonal automorphism $(x_1,...,x_n) \mapsto (x_1 \xi^{gq},x_2 \xi^{gm},x_3,...,x_n)=(x_1 \xi,x_2,...,x_n)$, which concludes the proof.
\proofend
\end{proof}

\begin{corollary}[Conjugacy of affine automorphisms of finite order]
{\it Two affine automorphisms of the same finite order are conjugate in $\BirAn$, for $n>1$.}
\proofend
\end{corollary}

\subsection{Conjugacy classes of diagonal automorphisms}
Let us now continue the work for diagonal automorphisms that are not necessary of finite order.

\begin{proposition}
\label{Prp:GroupActionDiag}
{\it 
Two diagonal automorphisms of $\An$ that are conjugate in $\BirAn$ are conjugate by an element of $\GLnZ$. }\\
\end{proposition}
\begin{proof}
To each diagonal automorphism $\rho(\alpha) \in \AutKn, \alpha=(\alpha_1,...,\alpha_n) \in (\K^{*})^n$, we associate the kernel $\KerL{\alpha}$ of the following homomorphism:
\begin{center}
$\begin{array}{cccc}
\delta_{\alpha}: & \mathbb{Z}^n & \rightarrow&\K^{*}\\
&(t_1,t_2,...,t_n)& \mapsto & {\alpha_1}^{t_1} {\alpha_2}^{t_2} ... {\alpha_n}^{t_n}.
\end{array}$\end{center}
It is easy to verify that $\delta_{M\bullet \alpha}=\delta_{\alpha} \circ \Mt$, for any $M\in \GLnZ$ and so that $\KerL{M\bullet \alpha}=\Mt^{-1}(\KerL{\alpha})$. If $\delta_{\alpha}$ is not injective, we can choose $M$ (by theorem on Smith's normal form) such that $\KerL{M\bullet \alpha}$ is generated by $k_1 e_1,k_2 e_2,...k_r e_r$, where $e_1,e_2,...,e_n$ are the canonical basis vectors of $\mathbb{Z}^n$, $r\leq n$ and the $k_i$ are positive integers such that $k_i$ divides $k_{i-1}$ for $i=2,...,r$. Let $\alpha'=(\alpha'_1,...,\alpha'_n)=M\bullet \alpha$. We observe that $\alpha'_i$ is a primitive $k_i$-th root of the unity for $i=1,...,r$. In particular, we have $\alpha'_2=(\alpha'_1)^s$ for an integer $s$, so $(s,-1,0,...0) \in \KerL{M\bullet \alpha}$, and then
 $k_2=k_3=..=k_r=1$. 

Let $\rho(\alpha)$ and $\rho(\beta)$ be two diagonal automorphisms conjugate in $\BirAn$.
We will assume that $\KerL{\alpha}$  is generated by $k e_1,e_2,... ,e_r$ (replacing $\rho(\alpha)$ by another element of its orbit if necessary, as we did above), so that $\alpha_1$ is a primitive $k$-th root of the unity and $\alpha_2=\alpha_3=...=\alpha_r=1$, if $r>0$. By the way, if $r=n$, we see another time (Proposition \ref{Prp:FiniteOrder}) that any diagonal automorphism of finite order is conjugate to another one of the form $(x_1,...,x_n) \mapsto (\xi x_1,x_2,...,x_n)$.

The condition of conjugation implies that there exists a birational map 
\begin{center}$\varphi:(x_1,x_2,...,x_n) \dasharrow (\varphi_1{(x_1,...,x_n)},\varphi_2{(x_1,...,x_n)},...,\varphi_n{(x_1,...,x_n)}),$\end{center}
 where $\varphi_i \in \K(X)$, such that $\varphi \circ \rho(\alpha)=\rho(\beta)  \circ \varphi$, and then 
\begin{center}
$\varphi_i{(\alpha_1 x_1,...,\alpha_n x_n)}=\rho(\alpha)^{*}(\varphi_i)=\beta_i \varphi_i$ for $i=1,...,n$.
\end{center}
So, every $\varphi_i$ must be an eigenvector of $\rho(\alpha)^{*}$, viewed as linear map of vector spaces, and since $\rho(\alpha)^{*}$ gives an automorphism of $\K[X]$, there must exist $F_i,G_i \in \K[X]$ eigenvectors of $\rho(\alpha)^{*}_ {|_{\K[X]}}$ such that $\varphi_i=\frac{F_i}{G_i}$. This map can be diagonalized along the basis of eigenvectors ${X_1}^{t_1}{X_2}^{t_2}...{X_n}^{t_n}$, and two vectors of this basis have the same eigenvalue if and only if the difference of the exponent vectors is in the kernel group $\Delta_{\alpha}$ associated to $\rho(\alpha)$.

By a choice of representants in the classes mod $\Delta_{\alpha}$, we can put the transformation $\varphi$ on the following form:
 \begin{center}
$\varphi:(x_1,...,x_n) \dasharrow ({x_1}^{a_{11}}\cdots{x_n}^{a_{1n}} \psi_1({x_1}^{k},x_2,...,{x_r}),...,{x_1}^{a_{n1}}\cdots {x_n}^{a_{nn}}\psi_n({x_1}^{k},x_2,...,{x_r}))$
\end{center}
where each $\psi_i$ belongs to $\K(X_1,...,X_r)$.

Hence, the map $\varphi$ defines a matrix $A=(a_{ij})_{i,j=1}^{n}$ in $Mat_{n,n}(\mathbb{Z})$, wich is not unique because of the choice of the representants. More precisely, $\varphi$ defines an unique element of the quotient of $Mat_{n.n}(\mathbb{Z})$ by the relation \begin{center}$(b_{ij})_{i,j=1}^{n} \sim  (c_{ij})_{i,j=1}^{n}$ if $b_{i1}-c_{i1} \in k \mathbb{Z}$  and $b_{ij}=c_{ij}$ for $i=1,...,n$ and $j>r$. \end{center}

We observe that $\beta=(\beta_1,...,\beta_n)=({{\alpha}_1}^{a_{11}}...{{\alpha}_n}^{a_{1n}},...,{{\alpha}_1}^{a_{n1}}...{{\alpha}_n}^{a_{nn}})$.
The birationaliy of $\varphi$ implies, (see Lemma 
\ref{Lem:birat-invertible} below), that there exists a matrix $B=(b_{ij})_{i,j=1}^{n}$ in the equivalence class associated to $\varphi$ which is invertible. So, we get from this that $T(B) \rho(\alpha) T(B)^{-1}=\rho(\beta)$, by the fact that 
\begin{center}$\begin{array}{cccccl}
 B\bullet \alpha&=(&{{\alpha}_1}^{b_{11}}\cdots{{\alpha}_n}^{b_{1n}}&,...,&{{\alpha}_1}^{b_{n1}}\cdots{{\alpha}_n}^{b_{nn}}&)\\
&=(&{{\alpha}_1}^{a_{11}}\cdots{{\alpha}_n}^{a_{1n}}&,...,&{{\alpha}_1}^{a_{n1}}\cdots{{\alpha}_n}^{a_{nn}}&)\\
&=(&\beta_1&,...,&\beta_n&)=\beta.\end{array}$\end{center} \proofend
\end{proof}

\begin{lemma}
\label{Lem:birat-invertible}
{\it If the rational map
 \begin{center}
$\begin{array}{rr}\varphi(x_1,...,x_n)&= (x_1^{a_{11}}\cdots x_n^{a_{1n}} \psi_1({x_1}^{k},x_2,...,{x_r}),...,\\
&x_1^{a_{n1}}\cdots x_n^{a_{nn}}\psi_n({x_1}^{k},x_2,...,{x_r}))\end{array}$
\end{center}
is birational, then there exists an invertible matrix $B=(b_{ij})_{i,j=1}^{n} \in \GLnZ$ such that
$b_{i1}-a_{i1} \in k \mathbb{Z}$ and $b_{ij}=a_{ij}$ for $i=1,...,n$ and $j>r$.}
\end{lemma}
\begin{proof}
First of all, if $n=1$ and $r=0$, the lemma is deduced from the fact that the map $x \dasharrow x^{a}$ is birational if and only if $a= \pm 1$.

Suppose that $r=n\geq1$ and let $\xi$ a primitive $k$-th root of the unity. The conjugation of the linear automorphism of finite order  $\nu_1:(x_1,...,x_n) \mapsto (x_1 \xi,x_2,...,x_n)$ by $\varphi$ gives the linear automorphism (of the same order) $\nu_2:(x_1,...,x_n) \mapsto (x_1 \xi_1^{a_{11}},...,x_n \xi^{a_{n1}})$. If $n=1$, the lemma follows from the fact that $x \mapsto \xi x$ and $x\mapsto \xi^{a_{11}} x$ are conjugated in $PGL(2,\K)$ if and only if $\xi^{a_{11}}=\xi^{\pm 1}$, so $a_{11}=\pm 1\ mod\ k$ . If $n>1$, by Proposition \ref{Prp:FiniteOrder}, there exists an invertible matrix $B=(b_{ij})_{i,j=1}^n \in \GLnZ$ such that $T(B)$ conjugates $\nu_1$ to $\nu_2$. Explicitly, it gives that $\xi^{b_{i1}}=\xi^{a_{i1}}$ for $i=1,...,n$, so $b_{i1}-a_{i1} \in k \mathbb{Z}$ for $i=1,...,n$. The lemma is then proved in this case.

It remains to prove the lemma when $r$ is strictly smaller than $n$ and $n>1$. In this case, the matrix $(a_{ij})_{1\leq i \leq n, r+1\leq j \leq n}$ has the maximal rank $n-r$, since $x_{r+1},...,x_n$ appear only in the monomial part.
There exists then an invertible matrix $M \in GL_n(\mathbb{Z})$ such that the coefficients of $A'=MA$ are the same as the ones of the identity, for columns $r+1,...,n$.
The composition of $T(M)$ and $\varphi$ gives a new birational map $T(M) \circ \varphi$:

 \begin{center}
$\begin{array}{rl}T(M)\circ \varphi(x_1,...,x_n)=(&x_1^{a_{11}'}\cdots x_{r}^{a_{1,r}'}\psi'_1(x_1^{k},x_2,...,x_r),...,\\
&x_1^{a_{r,1}'}\cdots x_{r}^{a_{r,r}'}\psi'_{r}({x_1}^{k},x_2,...,{x_r}),\\
&{x_1}^{a_{r+1,1}'}\cdots x_{r}^{a_{r+1,r}'}x_{r+1}\psi'_{r+1}({x_1}^{k},x_2,...,{x_r}),...,\\
& {x_1}^{a_{n,1}'}\cdots x_{r}^{a_{n,r}'}x_n \psi'_n({x_1}^{k},x_2,...,{x_r})),\end{array}$
\end{center}

which has the same structure as $\varphi$, with matrix $M A$.

 This map exchanges affine spaces of codimension $r$ of type $x_1=\tau_1,x_2=\tau_2,...,x_{r}=\tau_{r}$ and so must be its inverse. So, the map of $\K^{r}$ given by the $r$ first coordinates of $T(M)\circ \varphi$ must be birational.
 
By induction, there exists an invertible matrix $C'=(c'_{ij})_{i,j=1}^{n-1} \in GL_{r}(\mathbb{Z})$ such that $c'_{i1}-a'_{i1} \in k \mathbb{Z}$ for $i=1,...,r$. 
The matrix $B'=\left(\begin{array}{c|c} C' & 0  \\ \hline 0 & Id\end{array}\right)$ (where $Id \in GL_{n-r}(\mathbb{Z})$ denotes the identity matrix) satisfy the condition that $b'_{i1}-a'_{i1} \in k \mathbb{Z}$ and $b'_{ij}=a'_{ij}$ for $i=1,...,n$ and $j>r$. The same occurs for $B=M^{-1}B'$ and $A=M^{-1}A'$.
\proofend
\end{proof}

\subsection{Conjugacy classes of almost-diagonal automorphisms}
\label{subsection:Thm:GroupActionAlDiag}
The case of almost-diagonal automorphism is deduced from the situation of diagonal automorphisms:
\begin{proposition}
\label{Prp:Proof:Thm: GroupActionSemiDiag}
{\it Two almost-diagonal automorphisms of $\K^n$ are conjugate in the group $\BirAn$ if and only if they are conjugate by an element of $\GLnmZ$.}
\end{proposition}
\begin{proof}
Let $\theta(\alpha)$ and $\theta(\beta)$ be two almost-diagonal automorphisms of $\An$:
\begin{center}
$\theta(\alpha): (x_1,x_2,...,x_n) \mapsto (x_1+1, \alpha_2 x_2,\alpha_3 x_3,...,\alpha_{n} x_{n})$\\
$\theta(\beta): (x_1,x_2,...,x_n) \mapsto (x_1+1, \beta_2 x_2,\beta_3 x_3,...,\beta_{n} x_{n})$.
\end{center}
We suppose that  $\theta(\alpha)$ and $\theta(\beta)$ are conjugate in $\BirAn$, so that there exists a birational map \begin{center}
$\varphi:(x_1,x_2,...,x_n) \dasharrow (\frac{P_1(x_1,...,x_n)}{Q_1(x_1,...,x_n)},\frac{P_2(x_1,...,x_n)}{Q_2(x_1,...,x_n)},...,\frac{P_n(x_1,...,x_n)}{Q_n(x_1,...,x_n)})$
\end{center}, with $P_i,Q_i \in \K[X]$ without common divisor, such that $\varphi \circ \theta(\alpha)=\theta(\beta)  \circ \varphi$. This implies that 
\begin{center}
$\frac{P_i(x_1+1, \alpha_2 x_2,...,\alpha_{n} x_{n})}{Q_i(x_1+1, \alpha_2 x_2,...,\alpha_{n} x_{n})}=\beta_i \frac{P_i( x_1,...,x_n)}{Q_i( x_1,...,x_n)}$ for $i=2,...,n$.
\end{center}

Since $\theta(\alpha)$ induces an automorphism $\theta(\alpha)^{*}$  of $\K[X]$, all the $P_i$ and $Q_i$, for $i\geq 2$, must be eigenvectors of $\theta(\alpha)^{*}$, viewed as a $\K$-linear map. Such eigenvectors belongs to $\K[X_2,...,X_n]$ (Lemma $\ref{Lem:NoX1OnEigenVect}$), then the map 
$\varphi$ exchange lines of type $x_2=\tau_2,x_3=\tau_3,...,x_{n}=\tau_{n}$ and so must be its inverse. This implies that the map $\varphi':(x_2,x_3,...,x_n) \dasharrow (\frac{P_2(x_2,...,x_n)}{Q_2(x_2,...,x_n)},...,\frac{P_n(x_2,...,x_n)}{Q_n(x_2,...,x_n)})$, given by the $n-1$ last coordinates of $\varphi$ is birational.

Since the birational map $\varphi'$ of $\mathbb{K}^{n-1}$ conjugates the diagonal automorphisms
\begin{center}
$\rho(\alpha): (x_2,...,x_n) \mapsto (\alpha_2 x_2,\alpha_3 x_3,...,\alpha_{n} x_{n})$\\
$\rho(\beta)': (x_2,...,x_n) \mapsto (\beta_2 x_2,\beta_3 x_3,...,\beta_{n} x_{n})$,
\end{center}
Proposition $\ref{Prp:GroupActionDiag}$ shows that $\rho(\alpha)$ and $\rho(\beta)$ are conjugate by an element of $\GLnmZ \subset Bir(\K^{n-1})$. Then, $\theta(\alpha)$ and $\theta(\beta)$ are conjugated by the same element of $\GLnmZ\subset\BirAn$.
\proofend
\end{proof}

\section{Conjugacy classes of automorphisms of $\Pn$ in the Cremona group}
\label{Sec:ConjBirPn}
Let us work in $\Pn$, the projective $n$-space over $\K$.
Choosing a coordinate $x_i$, the open subset $U_i=\{(x_0:...:x_n) \in \Pn \ | \ x_i \not=0\}$ is isomorphic to $\An$ via the map $(x_0:...:x_n) \stackrel{\nu_i}{\mapsto} (\frac{x_0}{x_i},...,\frac{x_{i-1}}{x_i},\frac{x_{i+1}}{x_i},...,\frac{x_{n}}{x_i})$.
The restriction map $\varphi\mapsto{\nu_i}\varphi{\nu_i}^{-1}$ gives an isomorphism of the group $\BirPn$ of birational maps of $\Pn$ to $\BirAn$, that we call both the \defn{Cremona group}.

Let us recall that a birational map of $\Pn$ is given by a map $(x_0:...:x_n) \dasharrow (P_0(x_0,...,x_n):...:P_n(x_0,...,x_n))$, where $P_0,...,P_n \in \K[X_0,...,X_n]$ are homogeneous polynomials of the same degree (that will be called the degree of the map). As a birational map is biregular if and only if its degree is one, the group of automorphisms (biregular rational maps) of $\Pn$ is $\PGLnpl$. 

We will denote \textit{diagonal} (respectively \textit{almost-diagonal}) automorphisms of $\Pn$  maps of the form $(x_0:...:x_n) \mapsto (x_0:\alpha_1 x_1:...:\alpha_n x_n), \ \mbox{with } \alpha=(\alpha_1,...,\alpha_n) \in (\K^{*})^{n}$ (respectively $(x_0:...:x_n) \mapsto (x_0:x_0+x_1:\alpha_2 x_2:...:\alpha_n x_n), \ \mbox{with } \alpha=(\alpha_2,...,\alpha_n) \in (\K^{*})^{n-1}$).

It is clear that the restriction map  $\varphi\mapsto{\nu_0}\varphi{\nu_0}^{-1}$ gives an isomorphism of $\DKn$ (respectively $\ADKn$) on the group of diagonal (respectively the set of almost-diagonal) automorphisms of $\Pn$.

We can then explicit the conjugacy classes of automorphisms of $\Pn$ in the Cremona group using the work made in affine space:

Any automorphism of $\Pn$ is conjugate, in $\BirPn$, either to a diagonal automorphism
or to an almost-diagonal automorphism.
The conjugacy classes of diagonal and almost diagonal automorphisms of $\Pn$ in the Cremona group are given by the action of $\GLnZ$ on the diagonal automorphisms and of $\GLnmZ$ on the almost-diagonal automorphisms (see Section  \ref{Sec:ConjugationBirAn}).

Furthermore, if $n>1$,  two linear automorphisms of $\Pn$ of the same order are conjugate in $\BirPn$. This result was already proved in dimension $2$ in (\cite{bib:BBl}, Proposition 2.1.) with another method that works also in higher dimension.
\begin{acknowledgement}
  The author wishes to express his sincere gratitude to P. De La Harpe, I. Pan, F. Ronga and T.Vust for many interesting discussions during the preparation of this paper.  Moreover, the author acknowledges support from the {\it Swiss National Science Foundation}.
  %The author would like to thank the referee for his helpful criticisms and wishes to express his sincere gratitude to P. De La Harpe, I. Pan, F. Ronga and T.Vust for many interesting discussions during the preparation of this paper.  Moreover, the author acknowledges support from the {\it Swiss National Science Foundation}.
  \end{acknowledgement}

\bigskip

\end{document}